%
%

\documentstyle{amsppt}

\NoBlackBoxes

\pagewidth{5.0in}
\CenteredTagsOnSplits
\TagsOnRight

\define\mbox#1{\text{#1}}

\define\comp{{\Bbb C}}

\define\lie{{\Cal L}}
\define\sltc{{SL_2(\comp)}}

\define\gn#1{{\frak g_{#1}}}

\def\today{October 28, 1997}
\topmatter
\title  Frobenius Manifolds from Yang-Mills Instantons\\
\today
\endtitle
\rightheadtext{Frobenius Manifolds - \today}
\author Jan SEGERT \endauthor
\address Department of Mathematics,  
University of Missouri,
Columbia, MO 65211  \endaddress 
\email jan\@math.missouri.edu \endemail 

\abstract  We present an elementary self-contained account of 
semisimple Frobenius manifolds in three dimensions, 
and exhibit a new family of explicit examples.  These examples 
are constructed from Yang-Mills instantons with a certain symmetry.  
\endabstract

\endtopmatter
\document

\head 1.~Introduction \endhead

The concept of a Frobenius manifold was introduced and 
extensively developed by Dubrovin, whose lecture notes \cite{D} constitute 
the primary reference for Frobenius manifolds and many of the 
applications.   The lecture notes  of Hitchin \cite{Hi1} 
and of Manin \cite{Mn1} are also very good general references.  
Frobenius manifolds have appeared in a remarkably wide range of 
settings, including quantum cohomology \cite{RT},  
mirror symmetry and variation of Hodge structure \cite{G},  
unfoldings of singularities \cite{Au,Sa}, 
and the WDVV equation of topological quantum field theory  
 \cite{D}.  
Since Frobenius manifolds are relevant in the description of some  
deep geometrical phenomena, it is  not surprising that 
explicit solutions of the Frobenius manifold equations are 
rather difficult to construct.  

In this paper we present new Frobenius manifolds of dimension three. 
The paper is self-contained, no previous knowledge of Frobenius 
manifolds is assumed. 
In Section 2 we define 
semisimple Frobenius 
manifolds in the framework of canonical coordinates.  
In Section 3 we exhibit explicit formulae for the 
 new Frobenius manifolds.     
In Section 4 and the Appendix we use elementary Riemannian 
geometry to prove  some of the      
fundamental results for Frobenius manifolds \cite{D} for the special case 
of dimension three.  The cross-product effects numerous simplifications 
that are specific to dimension three.    
 Using  these  results, Proposition 2.2 in particular, 
it is easy to verify that the explicit formulae presented in Section 3 
are indeed Frobenius manifolds.       

The geometry used to construct the new Frobenius manifolds will 
be discussed in detail in another publication.  Here we  give only a 
brief outline.   
The theory of isomonodromic deformations offers an approach to 
the construction of semisimple Frobenius manifolds 
\cite{D,Hi,Mn1,Sa}.  
Isomonodromic deformations of a meromorphic connection on 
$\comp P^1$ are well-understood \cite{Ml}, 
but explicit examples are  difficult to construct.  
This difficulty is reflected in the correspondence between 
solutions of the Painlev\'e VI differential equation and a 
class of isomonodromic deformation \cite{F,JM};  Painlev\'e 
equations are notoriously difficult to solve explicitly.  
Hitchin \cite{Hi3,Hi2}   constructed some  
solutions of Painlev\'e VI by relating  certain isomonodromic deformations
to equivariant twistor geometry.  
The prototype is the  
irreducible linear $\sltc$ action on $\comp P^3$, 
which  gives rise to a natural 
flat meromorphic connection on $\comp P^3$.  The restriction to  
an embedded line $\comp P^1$ is isomonodromic under deformation of the line. 
The corresponding Frobenius manifold is the $n=0$ instance of Theorem 3.1. 

New isomonodromic transformations can in principle be 
constructed by applying ``Schlesinger transformations", which are meromorphic 
gauge transformations on $\comp P^1$ \cite{JM}.  
We refer to the recent paper  
of Manin \cite{Mn2} for a discussion of the corresponding 
transformations of  
Painlev\'e VI solutions.  
Equivariant twistor geometry provides a method for constructing 
some  Schlesinger transformations explicitly.  
The prototype  $\comp P^3$ with the irreducible $\sltc$ action 
is the equivariant twistor space of the Riemannian 
manifold $S^4$ with a certain isometric $SU_2$ action.  
The Atiyah-Ward correspondence \cite{At} relates 
anti-self-dual Yang-Mills instantons on $S^4$ to 
 certain holomorphic 
bundles on the twistor space $\comp P^3$.  An equivariant 
version  relates instantons with $SU_2$ symmetry to holomorphic 
bundles with $\sltc$ symmetry.  These equivariant objects   
were constructively classified 
in \cite{BS} by an equivariant version of the ADHM method \cite{ADHM}. 
Isomonodromic deformations can be generated from the equivariant ADHM 
monads, and all are related by   Schlesinger transformations. 
The Frobenius manifolds of Theorem 3.1 are constructed from 
the equivariant ADHM data.  

I would like to thank G. Bor and  N.J. Hitchin for helping me 
in understanding some of these topics.    

\head 2.~Semisimple Frobenius Manifolds\endhead

We define semisimple (or massive) Frobenius manifolds in the framework of  
local canonical coordinates \cite{D,Hi1,Mn1}, and state some basic 
results for three-dimensional Frobenius manifolds.  
All coordinates are complex, all functions are holomorphic, and 
all derivatives are  with respect to a complex variable.  
We follow the notational conventions of \cite{CDD; Sec. VI.A.4} for the exterior 
derivative $d$,  Lie derivative $\lie_V$, and interior product  $i_V$.  
Riemannian metrics are complex bilinear, not hermitian. 
 
We first introduce the Euler vector field and the identity vector field.  
The complex Lie group 
$\sltc$ acts on the Riemann sphere $\comp \, \cup \,  \{ \infty \}$ by 
fractional linear transformations.  The two-dimensional Borel subgroup  
$B \subset \sltc$ consisting  of the upper-triangular matrices is the 
stabilizer  of $\infty$.   $B$ acts by dilations and translations, 
$x \mapsto a \, x + b$, on 
the coordinates of a point $x \in \comp$.  
Using the notation $\partial = \frac \partial {\partial x}$, 
the dilations are generated by the 
``Euler vector field" $E = x \, \partial$ and the 
translations by the ``identity vector field" $I = \partial$, with $[I,E] = I$.  
More generally, for the action of 
 $B$ on  the coordinates $(x_1,x_2, \dots, x_n) \in \comp^n$ of 
an $n$-tuple of points in $\comp$,  the 
Euler vector field generating the dilations is the radial vector field 
$$\eqalign{ 
E &= x_1 \, \partial_1 + x_2 \, \partial_2 + \cdots + x_n \, \partial_n , \cr
}\tag 2.1
$$
and the identity vector field generating the translations  is 
$$\eqalign{ 
 I &= \partial_1 + \partial_2 + \cdots +\partial_n, \cr 
}\tag 2.2
$$
where $[I,E] = I$ as before.  Note that the Euler vector field is 
characterized by the property $\lie_E \, x_i = x_i$, and the 
identity vector field by the property $\lie_I \, x_i = 1$.  
A function $f$ on $\comp^n$, or on an 
open subset $U \subset \comp^n$,  will be called ``$B$-invariant" if 
$\lie_E \, f = 0$ and $\lie_I \, f = 0$.   
We will say that a function $f$ is of ``homogeneity $m$" if 
$\lie_E \, f = m \, f$ for a constant $m$. 
 
A Riemannian metric $g$ is ``flat" if the curvature of the associated 
Levi-Civita connection $\nabla$ is zero.   
A Riemannian metric on $U \subset \comp^n$ is ``diagonal" if 
it is of the form 
$$ g =   g_{11} \, dx_1 \otimes dx_1 + 
g_{22} \, dx_2 \otimes dx_2 + \cdots +g_{nn} \, dx_n \otimes dx_n.   \tag 2.3 $$

\proclaim{Definition 2.1} A ``semisimple Frobenius manifold" 
structure of homogeneity $m$ on 
open subset $U \subset \comp^n$  with 
``canoncial coordinates" $(x_1,x_2,\dots,x_n)$ consists of a diagonal 
metric $g$  satisfing the three conditions:    
\roster
\item"(M1)" 
$g$ is flat.     
\item"(M2)" 
The components $g_{ii}$ are functions of homogeneity $m$.  
\item"(M3)"
The identity vector field $I$ is covariantly constant with respect to 
the Levi-Civita connection.   
\endroster
\endproclaim

\noindent   
A Frobenius metric is ``nontrivial" if the components 
$g_{ii}$ are not all constant; the standard Euclidean metric 
is an example of a trivial Frobenius metric.        
Now apply definition 2.1 to an atlas of local 
coordinate charts on a manifold.    
A complex manifold $M$ with a complex metric $g$ 
and vector fields $E$ and $I$ satisfying $[I,E] = I$ 
is a ``semisimple Frobenius manifold" of homogeneity $m$ 
if every point of $M$ has a neighborhood 
which admits local canonical coordinates $(x_1,x_2,..,x_n)$ as above. 
It is convenient to relax this definition, requiring only 
that every point on some dense open subset has such a 
neighborhood, and that $g$ is nondegenerate and nonsingular only on 
some dense open subset.     

We now focus on  the three-dimensional case.      
The cross-ratio  
$$ t = { x_3 - x_1 \over x_2 - x_1 }  \tag 2.4 $$
is a $B$-invariant function on $\comp^3$.  In section 4, we will prove:  

\proclaim{Proposition 2.2} 
The metric $ g = g_{11} \, dx_1 \otimes dx_1 + 
g_{22} \, dx_2 \otimes dx_2 + g_{33} \, dx_3 \otimes dx_3 $ 
is a  homogeneity-$0$ Frobenius metric if and only if 
$$ 
t \, d\, g_{11}  =  
(1-t)\, d \,g_{22} = 
t\, (t-1)\, d \, g_{33} = 
2 \, c \, \sqrt{ g_{11} \, g_{22} \, g_{33}} \, dt  
\tag 2.5
$$
for some constant $c$.  The Frobenius metric is nontrivial 
if and only if $c$ is nonzero.   
\endproclaim
\noindent It is evident that such a $c$ is unique, and that the ``trace" 
$$ k = - {\frac {c^2} 2} \, (g_{11} + g_{22} + g_{33})  \tag 2.6
$$
of homogeneity-$0$ Frobenius metric $g$ is a 
constant.  The trace is  
unchanged under the rescaling $g \mapsto \alpha \, g$ by a constant 
$\alpha$.     

The two-dimensional Lie group $B$ is a symmetry group of 
three-dimensional Frobenius metrics.    
For a nontrivial homogeneity-$0$ Frobenius metric, 
$i_Y$ applied to eq.(2.6)  yields  
$$ 
t \, \lie_Y \, g_{11}  =  
(1-t)\, \lie_Y\,g_{22} = 
t\, (t-1)\, \lie_Y \, g_{33} = 
2 \, c \, \sqrt{ g_{11} \, g_{22} \, g_{33}} \, \lie_Y \, t  
$$
for any vector field $Y$.  
The $B$-invariance of $t$, $\lie_E \, t =\lie_I \, t =0$, then implies
the $B$-invariance of $g_{ii}$, 
$\lie_E \, g_{ii}=\lie_I \, g_{ii} =0$.  More generally for 
homogeneity $m$, condition (M2) states  
$\lie_E \, g_{ii} = m \, g_{ii}$, and the results of section 4 give 
$\lie_I \, g_{ii} = 0$.   
In the language of Riemannian geometry, $I$ is a Killing vector, 
$\lie_I \, g = 0$, and   $E$ is a weight-$(m+2)$ conformal Killing vector, 
$\lie_E \, g = (m+2) \, g$.

\head 3.~Frobenius Manifolds Constructed from Instantons \endhead

We start by constructing an atlas of canonical coordinate charts, 
and the corresponding Euler and identity 
vector fields, on a 
certain  hypersurface $M \subset \comp^4$.  
Let $(z_1,z_2,z_3,r)$ be the 
linear coordinates on $\comp^4$,  let $M$  
be the hypersurface defined by the vanishing of 
$$ Q = 
  {{\left( {z_1} - {z_2} \right) }^2} + 
  {{\left( {z_2} - {z_3} \right) }^2} + {{\left( {z_3}-{z_1}   \right) }^2}
-2\,{r^2},
$$
and let $j: M \to \comp^4$ denote the inclusion map.  
Every point in an open dense subset of $M$ 
has a neighborhood on which the  
restrictions $(x_1,x_2,x_3) = (j^*q_1, j^*q_2, j^*q_3)$ of 
the three functions  
$$\eqalign{q_1 &= 
2\,r\,\left( 2\,{z_1} - {z_2} - {z_3} \right) -3\,{{{z_1}}^2} - 6\,{z_2}\,{z_3},
\cr
q_2 &= 
2\,r\,\left( 2\,{z_2} - {z_3} - {z_1} \right) -3\,{{{z_2}}^2} - 6\,{z_3}\,{z_1},
\cr
q_3 &= 
2\,r\,\left( 2\,{z_3} - {z_1} - {z_2} \right) -3\,{{{z_3}}^2} - 6\,{z_1}\,{z_2}
\cr
}$$
define  local coordinates; the points which fail to have this 
property are 
characterized by the vanishing of the  Jacobian determinant.  
Of course 
these coordinates are only valid locally because the
 functions $x_i: M \to \comp$ are not one-to-one.    
The group of dilations on $\comp^4$ is generated by the 
radial vector field $2 \,\tilde E$, where 
$$ \eqalign{\tilde E &=\frac 1 2 \left( z_1 \, \frac \partial {\partial z_1} + 
z_2 \, \frac \partial {\partial z_2} + 
z_3 \,\frac \partial {\partial z_3}+
r \,\frac \partial {\partial r}
\right).  \cr
}
$$ 
Now $\tilde E$ is tangent to $M$ since $\lie_{\tilde E} \, Q = Q$, so 
the restriction of $\tilde E$ projects to a vector field $E$ on $M$. 
The vector field 
$E$ is the Euler vector field relative to  the local coordinates 
 since  
$\lie_{\tilde E} \, q_i = q_i$.    
The vector field 
$$ \eqalign{
\tilde I &= {-1 \over 6\, ( z_1 + z_2 + z_3) }\, \left(
\frac \partial {\partial z_1} + 
\frac \partial {\partial z_2} + 
\frac \partial {\partial z_3}\right)  \cr
}
$$ 
is tangent to $M$ since 
$\lie_{\tilde I} \, Q = 0$, so the restriction of $\tilde I$ projects  
to a vector field $I$ on $M$.   
The vector field 
$I$ is the identity vector field relative to  the local coordinates 
since 
$\lie_{\tilde I} \, q_i = 1$, which follows from  the identity 
$q_1 + q_2 + q_3 = -3 \, (z_1 + z_2 + z_3)^2$.  
Let the symmetric group $S_3$  act on 
$\comp^4$ by permuting the first three coordinates 
$(z_1,z_2,z_3)$. 
The  polynomial $Q$ is $S_3$-invariant, so the action maps 
the hypersurface $M$ to itself.  The $S_3$ action permutes the 
three functions $q_i$, leaving the vector fields $\tilde E$ and $\tilde I$ 
invariant.

Our main result  is the existence of a family of Frobenius metrics on $M$, 
the coefficients of which are   rational functions on 
$\comp^4$:          
  
\proclaim{Theorem 3.1} For each nonnegative integer $n$, there exist  
 triplets $( b_1, b_2, b_3 )$ and  $( u_1, u_2, u_3 )$ of 
explicitly computable homogeneous polynomials  
of degree $l \le 2(n^2 + n + 2)$  on $\comp^4$ such that  
$$  g = j^* \left( 
{ u_1 \over b_1} \, dq_1 \otimes dq_1
+
{ u_2 \over b_2} \, dq_2 \otimes dq_2
+
{ u_3 \over b_3} \, dq_3 \otimes dq_3
\right)  \tag 3.1 
$$
is the metric of homogeneity-$0$ Frobenius manifold on the 
hypersurface $M \subset \comp^4$.  The Frobenius metric $g$ 
is nontrivial, and has trace  
$k = \frac 1 2 \, (n + \frac 1 2)^2$.   The symmetric group 
$S_3$ acts by permutation on each of the 
triplets $(b_1,b_2,b_3)$, $(u_1,u_2,u_3)$, and $(q_1,q_2,q_3)$, so  
the Frobenius structure on $M$ is $S_3$-invariant. 
\endproclaim 

\noindent 
The constructive geometric proof of theorem 3.1 for 
all nonnegative integers $n$, using the  classification 
by Bor and the author \cite{BS} of  Yang-Mills instantons with a 
certain $SU_2$-symmetry, will be described elsewhere.  In the 
present paper, we exhibit  
the  polynomials $b_i$ and $u_i$  for $n \le 2$, and compute the
canonical coordinate expressions of the Frobenius metrics.   
Applying Proposition 2.2 to these expressions constitutes a computational  
proof of theorem 3.1 for $n \le 2$.    
We do not continue beyond $n=2$ because the size of the expressions grows 
very quickly with $n$.

The local canonical coordinate expressions of the form eq.(2.3) are easily 
evaluated for the Frobenius metrics eq.(3.1).   We first observe that      
a $B$-invariant function on $U \subset \comp^3$ depends only   
 on the cross-ratio $t$.  
However, a $B$-invariant function on $M$ is a  possibly 
{\it multi-valued} function of $t$, because the cross ratio 
$t:M \to \comp$, although well-defined,  is not one-to-one.  
If a path    
$ \gamma: \comp \to M$ is transverse to $B$-orbits, then it follows 
from Theorem 3.1 and $B$-invariance that  
$$ g_{ii}(x_1,x_2,x_3) = { u_i ( \gamma(w)) \over b_i (\gamma(w))} $$
where $w$ is a (possibly non-unique) solution of 
$$ {x_3 - x_1 \over x_2 - x_1} = 
{ q_3(\gamma(w)) - q_1(\gamma(w)) \over q_2(\gamma(w)) - q_1(\gamma(w))}. 
$$
Choosing the polynomial path    
$$  \gamma(w) =
 \left( w^2 -1  \, , \,  -2 \, w + 2 \, ,
  \, 2 \, w + 2 \,  , \,{w^2}+3 \right) \in M \subset \comp^4, $$
the $g_{ii}$ become rational functions of  degree at 
most $2 \, l$ in the variable $w$, where $w$ is a solution of  
$$
\eqalign{  {x_3 - x_1 \over x_2 - x_1}
 &=
{\left( w+1 \right)\,\left( w-3 \right) ^3 \over 
\left( w-1 \right) \,\left( w+3 \right)^3 }. \cr 
} $$
In terms of the cross-ratio $t$ of the local coordinates $(x_1,x_2,x_3)$, 
$w$ is a root of the following quartic polynomial 
with coefficients  depending on $t$: 
$$
\eqalign{  
\left( w-1 \right) \,\left( w+3 \right)^3 \,t - 
\left( w+1 \right)\,\left( w-3 \right) ^3 &= 0. \cr 
} \tag 3.2$$
We recall that the 
roots of a quartic polynomial can be expressed as explicit (multivalued)  
algebraic function of the coefficients by a formula analogous 
to the familiar quadratic formula, albeit much more 
complicated.  

We now exhibit the data of theorem 3.1 for the first few values of $n$.  
For $n =0$, the data is constructed from the  ADHM monad of the trivial    
Yang-Mills instanton.  The trivial instanton has instanton number $0$, 
and has $SU_2$ symmetry of \cite{BS}.   
The homogeneous polynomials $u_i$ and $b_i$ on $\comp^4$ have 
degree $l=2$:  
$$\eqalign{{ b_1} =  
    36\,\left( {z_3}-{z_1}  \right) \, \left( {z_1} - {z_2} \right),\qquad
{ u_1 } &=  \left( r  - ({z_3}-{z_1})+({z_1} - {z_2})  \right)^2
, \cr
{b_2} =  
    36\,\left( {z_1}-{z_2}  \right) \, \left( {z_2} - {z_3} \right),
\qquad
{ u_2 } &=  \left( r  - ({z_1}-{z_2})+({z_2} - {z_3})  \right)^2,
  \cr
{  b_3} = 
    36\,\left( {z_2}-{z_3}  \right) \, \left( {z_3} - {z_1} \right),
\qquad
{ u_3 } &=  \left( r  - ({z_2}-{z_3})+({z_3} - {z_1})  \right)^2.
\cr
} \tag 3.3
$$
The $S_3$ symmetry of the 
triplet $(b_1,b_2,b_3)$ and of the triplet $(u_1,u_2,u_3)$ is evident.  
The canonical coordinate expression of the $n=0$ metric is 
$$\eqalign{ g &= 
 -{\frac{ \left( w-1 \right) \,\left( w+1 \right)   }
    {4\,\left( w-3 \right) \,\left( w+3 \right) }}\, dx_1 \otimes dx_1 
  -{\frac{\left( w-1 \right) }{4\,w\,\left( w+3 \right) }}\, 
dx_2 \otimes dx_2 \cr
&\quad +
  {\frac{ \left(w+1\right)}{4\,w\,\left( w-3 \right) }} \,dx_3 \otimes dx_3, \cr
} 
\tag 3.4
$$
where $w$ is related to the cross-ratio $t$ by eq.(3.2).  
The reader may easily check that this metric satisfies the conditions of 
Proposition 2.2 with $c=1$, and has trace $k=\frac 1 8$.

For $n =1$, the data is constructed from the  ADHM monad of the 
``basic instanton" \cite{At}, which is the 
unique $SU_2$-symmetric instanton with instanton number $1$.   
The homogeneous polynomials $u_i$ and $b_i$ on $\comp^4$ have 
degree $l=6$:
$$ \eqalign{b_1 &=
 {216\, {r^2}\, {{\left( {z_2} - {z_3} \right) }^2} 
\,\left( {z_3} -{z_1}  \right)\,\left( {z_1} - {z_2} \right)  },\cr}
$$
and  $u_1$ is an irreducible   
polynomial with 84 terms:  
$$ \eqalign{
u_1 & =  
-3\,{r^6} + 14\,{r^5}\,{z_1} - 21\,{r^4}\,{{{z_1}}^2} + 
     4\,{r^3}\,{{{z_1}}^3} + 19\,{r^2}\,{{{z_1}}^4} - 18\,r\,{{{z_1}}^5} + 
     5\,{{{z_1}}^6} - 7\,{r^5}\,{z_2} \cr
&+ 21\,{r^4}\,{z_1}\,{z_2} - 
     6\,{r^3}\,{{{z_1}}^2}\,{z_2} - 38\,{r^2}\,{{{z_1}}^3}\,{z_2} + 
     45\,r\,{{{z_1}}^4}\,{z_2} - 15\,{{{z_1}}^5}\,{z_2} + 
     39\,{r^4}\,{{{z_2}}^2} \cr
&+ 108\,{r^3}\,{z_1}\,{{{z_2}}^2} + 
     240\,{r^2}\,{{{z_1}}^2}\,{{{z_2}}^2} + 
     66\,r\,{{{z_1}}^3}\,{{{z_2}}^2} + 33\,{{{z_1}}^4}\,{{{z_2}}^2} - 
     53\,{r^3}\,{{{z_2}}^3} \cr
&- 221\,{r^2}\,{z_1}\,{{{z_2}}^3} - 
     144\,r\,{{{z_1}}^2}\,{{{z_2}}^3} - 41\,{{{z_1}}^3}\,{{{z_2}}^3} + 
     22\,{r^2}\,{{{z_2}}^4} + 27\,r\,{z_1}\,{{{z_2}}^4} - 
     33\,{{{z_1}}^2}\,{{{z_2}}^4} \cr
&+ 12\,r\,{{{z_2}}^5} + 
     51\,{z_1}\,{{{z_2}}^5} - 10\,{{{z_2}}^6} - 7\,{r^5}\,{z_3} + 
     21\,{r^4}\,{z_1}\,{z_3} - 6\,{r^3}\,{{{z_1}}^2}\,{z_3} - 
     38\,{r^2}\,{{{z_1}}^3}\,{z_3} \cr
&+ 45\,r\,{{{z_1}}^4}\,{z_3} - 
     15\,{{{z_1}}^5}\,{z_3} - 99\,{r^4}\,{z_2}\,{z_3} - 
     204\,{r^3}\,{z_1}\,{z_2}\,{z_3} - 
     366\,{r^2}\,{{{z_1}}^2}\,{z_2}\,{z_3} \cr
&- 312\,r\,{{{z_1}}^3}\,{z_2}\,{z_3} + 9\,{{{z_1}}^4}\,{z_2}\,{z_3} + 
     51\,{r^3}\,{{{z_2}}^2}\,{z_3} + 
     183\,{r^2}\,{z_1}\,{{{z_2}}^2}\,{z_3} \cr
&+ 234\,r\,{{{z_1}}^2}\,{{{z_2}}^2}\,{z_3} - 
     9\,{{{z_1}}^3}\,{{{z_2}}^2}\,{z_3} + 133\,{r^2}\,{{{z_2}}^3}\,{z_3} + 
     180\,r\,{z_1}\,{{{z_2}}^3}\,{z_3} + 
     255\,{{{z_1}}^2}\,{{{z_2}}^3}\,{z_3} \cr
&- 87\,r\,{{{z_2}}^4}\,{z_3} - 
     189\,{z_1}\,{{{z_2}}^4}\,{z_3} + 9\,{{{z_2}}^5}\,{z_3} + 
     39\,{r^4}\,{{{z_3}}^2} + 108\,{r^3}\,{z_1}\,{{{z_3}}^2} + 
     240\,{r^2}\,{{{z_1}}^2}\,{{{z_3}}^2} \cr
&+  66\,r\,{{{z_1}}^3}\,{{{z_3}}^2} + 33\,{{{z_1}}^4}\,{{{z_3}}^2} + 
     51\,{r^3}\,{z_2}\,{{{z_3}}^2} + 
     183\,{r^2}\,{z_1}\,{z_2}\,{{{z_3}}^2} + 
     234\,r\,{{{z_1}}^2}\,{z_2}\,{{{z_3}}^2}\cr
& -  9\,{{{z_1}}^3}\,{z_2}\,{{{z_3}}^2} - 
     291\,{r^2}\,{{{z_2}}^2}\,{{{z_3}}^2} - 
     504\,r\,{z_1}\,{{{z_2}}^2}\,{{{z_3}}^2} - 
     369\,{{{z_1}}^2}\,{{{z_2}}^2}\,{{{z_3}}^2} + 
     84\,r\,{{{z_2}}^3}\,{{{z_3}}^2}\cr
& + 123\,{z_1}\,{{{z_2}}^3}\,{{{z_3}}^2} + 72\,{{{z_2}}^4}\,{{{z_3}}^2} - 
     53\,{r^3}\,{{{z_3}}^3} - 221\,{r^2}\,{z_1}\,{{{z_3}}^3} - 
     144\,r\,{{{z_1}}^2}\,{{{z_3}}^3} - 41\,{{{z_1}}^3}\,{{{z_3}}^3} \cr
&+  133\,{r^2}\,{z_2}\,{{{z_3}}^3} + 180\,r\,{z_1}\,{z_2}\,{{{z_3}}^3} + 
     255\,{{{z_1}}^2}\,{z_2}\,{{{z_3}}^3} + 
     84\,r\,{{{z_2}}^2}\,{{{z_3}}^3} + 
     123\,{z_1}\,{{{z_2}}^2}\,{{{z_3}}^3}\cr
& - 137\,{{{z_2}}^3}\,{{{z_3}}^3} + 
     22\,{r^2}\,{{{z_3}}^4} + 27\,r\,{z_1}\,{{{z_3}}^4} - 
     33\,{{{z_1}}^2}\,{{{z_3}}^4} - 87\,r\,{z_2}\,{{{z_3}}^4} - 
     189\,{z_1}\,{z_2}\,{{{z_3}}^4} \cr
&+ 72\,{{{z_2}}^2}\,{{{z_3}}^4} + 
     12\,r\,{{{z_3}}^5} + 51\,{z_1}\,{{{z_3}}^5} + 9\,{z_2}\,{{{z_3}}^5} - 
     10\,{{{z_3}}^6}. \cr
}
$$
The canonical coordinate expression of the $n=1$ metric is 
$$\eqalign{
g &= -{\frac{9\,{{\left( w-1 \right) }^3}\,{{\left( w+1\right) }^3}}
    {4\,\left( w-3 \right) \,\left( w+3 \right) \,
      {{\left(  {w^2} +3\right) }^2}}} \, dx_1\otimes dx_1
  -{\frac{9\,{{\left( w-1 \right) }^3}}
    {w\,\left( w+3 \right) \,{{\left( {w^2}+3 \right) }^2}}}\, 
dx_2 \otimes dx_2 \cr
& \quad +
  {\frac{9\,{{\left( w+1 \right) }^3}}
    {w\,\left( w-3\right) \,{{\left(  {w^2} +3\right) }^2}}} \,dx_3\otimes dx_3. 
\cr
}
$$
This metric satisfies the conditions of Proposition 2.2 with $c=1$, and 
has trace $k = \frac 9 8$.    

For $n =2$, the data is constructed from the  ADHM monad of 
the unique $SU_2$-symmetric instanton \cite{BS} with instanton number $3$.   
The homogeneous polynomials $u_i$ and $b_i$ on $\comp^4$ have 
degree $l=10$:
$$\eqalign{ b_1 &= 
36\,\left( {z_2} - {z_3} \right) ^2 
\,\left( {z_3} -{z_1}  \right)\,\left( {z_1} - {z_2} \right)  \, 
 ( 52\,{r^3} + 2\,{{{z_1}}^3} - 
       3\,{{{z_1}}^2}\,{z_2} - 3\,{z_1}\,{{{z_2}}^2} \cr
&\quad +  2\,{{{z_2}}^3} - 3\,{{{z_1}}^2}\,{z_3} + 
       12\,{z_1}\,{z_2}\,{z_3} - 3\,{{{z_2}}^2}\,{z_3} - 
       3\,{z_1}\,{{{z_3}}^2} - 3\,{z_2}\,{{{z_3}}^2} + 
       2\,{{{z_3}}^3} )^2, \cr
}
$$
and $u_1$ is an irreducible polynomial with 
283 terms, of which we only exhibit the first few: 
$$
\eqalign{u_1 &= 
25\,( 174\,{r^{10}} - 764\,{r^9}\,{z_1} + 1122\,{r^8}\,{{{z_1}}^2} - 
    244\,{r^7}\,{{{z_1}}^3} - 998\,{r^6}\,{{{z_1}}^4} + 1008\,{r^5}\,{{{z_1}}^5} 
\cr
&-  274\,{r^4}\,{{{z_1}}^6} - 20\,{r^3}\,{{{z_1}}^7} - 24\,{r^2}\,{{{z_1}}^8} + 
    20\,r\,{{{z_1}}^9} + 382\,{r^9}\,{z_2} - 1122\,{r^8}\,{z_1}\,{z_2} 
\cr
&+ 366\,{r^7}\,{{{z_1}}^2}\,{z_2} + 1996\,{r^6}\,{{{z_1}}^3}\,{z_2} - 
    2520\,{r^5}\,{{{z_1}}^4}\,{z_2} + 822\,{r^4}\,{{{z_1}}^5}\,{z_2} + 
\cdots \cdots) .
\cr
}
$$
The canonical coordinate expression of the $n=2$ metric is 
$$\eqalign{
g &= -{\frac{25\,{{\left( w-1 \right) }^5}\,{{\left( w+1 \right) }^5}\,
      {{\left( {w^2}+5 \right) }^2}}{4\,\left( w-3 \right) \,
      \left( w+3 \right) \,{{\left( {w^2} +1\right) }^2}\,
      {{\left(  {w^2} - 2\,w + 5\right) }^2}\,
      {{\left(  {w^2} + 2\,w + 5  \right) }^2}}}\, dx_1 \otimes dx_1 \cr
&\quad -
  {\frac{25\,{{\left( w-1 \right) }^5}\,
      {{\left( 3\,{w^2} + 2\,w + 7 \right) }^2}}{w\,\left( w+3 \right) \,
      {{\left( {w^2} +1\right) }^2}\,
      {{\left(  {w^2}- 2\,w + 5 \right) }^2}\,
      {{\left( {w^2}  + 2\,w + 5 \right) }^2}}}\, dx_2 \otimes dx_2 \cr
&\quad +
  {\frac{25\,{{\left( w+1 \right) }^5}\,
      {{\left( 3\,{w^2} - 2\,w + 7 \right) }^2}}{w\,\left( w-3 \right) \,
      {{\left( {w^2} +1\right) }^2}\,
      {{\left( {w^2} - 2\,w + 5 \right) }^2}\,
      {{\left( {w^2} + 2\,w + 5 \right) }^2}}} \, dx_3 \otimes dx_3. \cr
}
$$
This metric satisfies the conditions of Proposition 2.2 with $c=1$, 
and has trace $k = \frac {25} 8$.    
For reasons of brevity, we do not continue beyond $n=2$.

Corollary 4.3 below  associates a pair of 
Frobenius metrics $g^\pm$ of homogeneity $m = \pm \sqrt{8 \, k}$ to each  
Frobenius metric  $g$ of homogeneity $0$ and trace $k$.  
The metrics $g^\pm$ associated to the $n=0$ 
metric of Theorem 3.1 can be expressed in the form eq.(3.1).  
The homogeneous polynomials $u_i^+$ have degree $4$,    
$$ u_1^+ = 
{{{\left( {r^2} + 4\,r\,{z_1} - 5\,{{{z_1}}^2} - 2\,r\,{z_2} + 
         5\,{z_1}\,{z_2} + {{{z_2}}^2} - 2\,r\,{z_3} + 5\,{z_1}\,{z_3} - 
         7\,{z_2}\,{z_3} + {{{z_3}}^2} \right) }^2}} , 
$$
the polynomials $u_i^-$ have degree $0$, $u_1 = 1$, and  
the $b_i^\pm$ are equal to the degree $2$ polynomials $b_i$ of eq.(3.4).  
The canonical coordinate expressions are  
$$\eqalign{g^+ &= 
 -{\frac{4\,\left( w+1 \right)  }
    {\left( w-3  \right) \,{{\left( w+3 \right) }^4}}} 
\,\left( {x_1} - {x_2} \right)\, dx_1 \otimes dx_1
  -{\frac{ {{\left( w+1 \right) }^4}}
    {4\,w\,{{\left( w+3 \right) }^4}}} 
\, \left( {x_1} - {x_2} \right) \,dx_2 \otimes dx_2 \cr
&\quad + {\frac{{{\left( w-1 \right) }^3}
   \,\left( w+1 \right) 
       }{4\,w\,\left( w-3 \right) \,
      {{\left( w+3 \right) }^3}}} 
\,\left( {x_1} - {x_2} \right)\, dx_3 \otimes dx_3, \cr
g^- &= 
 -{\frac{{{\left( w+3 \right) }^2}}
     {\left( w-3 \right) 
\,\left( w+1 \right) }} \, \left( {x_1} - {x_2} \right)^{-1}\, dx_1 \otimes dx_1 
  -{\frac{{{\left( w+3 \right) }^2}}
    {4\,w }} \, \left( {x_1} - {x_2} \right)^{-1}\,dx_2 \otimes dx_2 
\cr 
&\quad +
  {\frac{\left( w-1 \right) \,{{\left( w+3 \right) }^3}}
    {4\, w\,\left( w-3 \right) 
\,\left( w+1 \right)  }} 
\, \left( {x_1} - {x_2} \right)^{-1}\, dx_3 \otimes dx_3 .\cr
}
$$
The canonical coordinate expressions for the 
Frobenius metrics of homogeneity 
$m= \pm 3, \pm 5, \dots$ associated to $n = 1,2, \dots$ can also be 
computed explicitly, but I do not know whether these metrics 
can be expressed in the form eq.(3.1).

Applying the 
correspondence between semisimple three-dimensional 
Frobenius manifolds  and Painlev\'e transcendents \cite{D,Hi1,JM} 
to the Frobenius manifolds 
of theorem 3.1, we obtain for 
each nonnegative integer $n$  
 two distinct solutions $\lambda_\pm(t)$ of the Painlev\'e VI equation 
$$ \eqalign{
{d^2 \lambda_\pm \over dt^2} = &
{1 \over 2} \left({1 \over \lambda_{\pm}} + {1 \over \lambda_\pm -1} + {1 \over 
\lambda_\pm-t}\right)
\left( {d\lambda_\pm \over dt} \right)^2 
- \left( {1 \over t} + {1 \over t -1} + {1 \over \lambda_\pm - t} \right) 
{d\lambda_\pm \over dt} \cr
&+ { \lambda_\pm(\lambda_\pm-1)(\lambda_\pm-t) \over t^2 (t-1)^2} 
\left(\alpha_\pm + \beta{t \over \lambda_\pm^2} + \gamma{t-1\over(\lambda_\pm-1)^2} + 
\delta{t (t-1) \over (\lambda_\pm-t)^2}\right) , \cr} 
$$
with  
$$\eqalign{   
\alpha_\pm & = 
{\textstyle{1 \over 2}}{\left( (n+{\textstyle{1 \over 2}}) \mp 1 \right)^2 }, \quad
\beta = -{\textstyle{1 \over 2}}{\left( n+{\textstyle{1 \over 2}}\right)^2},\quad 
\gamma = {\textstyle{1 \over 2}}{\left( n+{\textstyle{1 \over 2}}\right)^2 },\quad
\delta  ={\textstyle{1 \over 2}} -
{\textstyle{1 \over 2}}{\left( n+{\textstyle{1 \over 2}}\right)^2 }. 
\cr
}
$$
The solution $\lambda_\pm$ is a rational function of degree
$l_\pm \le  2( n^2 + n + 2)$ in the variable $w$, where 
$w$ is related to $t$ as in eq.(3.3) above. 
The explicit formulae for $\lambda_\pm$ are 
 exhibited in \cite{Se} for $n \le 4$, building on  
Hitchin's previous computation \cite{Hi3} of $\lambda_+$ for $n=0$.

\head 4.~Frobenius Coframes in Dimension Three\endhead

An ``orthonormal frame" for a Riemannian metric $g$ is a triplet 
of vector fields $e_i$ such that 
$ g(e_i, e_j) = \delta_{i,j} $. 
The dual ``orthonormal coframe" is the triplet  
of one-forms $\theta_j$ such that 
$ \langle \theta_i , e_j \rangle = \delta_{i,j} $. 
The metric can be reconstructed from an orthonormal coframe by  
$ g = \theta_1 \otimes \theta_1 + \theta_2 \otimes \theta_2 
+ \theta_3 \otimes \theta_3 $.  We will follow the notational conventions  
of \cite{CDD; Sec. VI.B}.  
 
A linear connection $\nabla$ on the tangent bundle which preserves  
the metric is an ``orthogonal" connection.      
Relative to an orthonormal frame, an 
orthogonal connection $\nabla$ is 
expressed in terms of the triplet of connection one-forms $\Omega_k$ as follows: 
If $X = X_1 \, e_1 + X_2 \, e_2 + X_3 \, e_3$, then 
$\nabla X = (\nabla X)_1 \otimes e_1 + (\nabla X)_2 \otimes e_2 + 
(\nabla X)_3 \otimes e_3$, where 
$$ \eqalign{
(\nabla X)_1 &= d X_1 - \Omega_2 \, X_3 + \Omega_3 \, X_2 \cr
 (\nabla X)_2  &= d X_2 - \Omega_3 \, X_1 + \Omega_1 \, X_3 \cr
(\nabla X)_3 &=  d X_3 - \Omega_1 \, X_2 + \Omega_2 \, X_1. \cr
} \tag 4.1
$$
using the three-vector notation   
$$
\vec X = \pmatrix 
X_1 \cr X_2 \cr X_3 \cr \endpmatrix, \qquad
\vec \Omega = \pmatrix \Omega_1 \cr \Omega_2 \cr \Omega_3 \cr 
\endpmatrix, \qquad 
 \vec \theta = \pmatrix \theta_1 \cr \theta_2 \cr \theta_3 \cr 
\endpmatrix,
$$
and the cross-product $\times$, eq.(4.1) becomes
$$ \vec \nabla \vec X = d \, \vec X - \vec \Omega \times \vec X . $$
The Levi-Civita connection of $g$ is the unique 
orthogonal connection with zero torsion.  
If $\vec \theta$ is an orthonormal coframe for $g$, then 
the connection form $\vec \Omega$ corresponding to the 
Levi-Civita connection is the unique solution of the 
Cartan torsion equation 
$$ d \, \vec \theta - \vec \Omega \times \vec \theta = 0 .$$
We will say that $\vec \Omega$ is the ``Levi-Civita connection form"  
of the orthononormal coframe $\vec \theta$ if the Cartan 
torsion equation holds.   
The curvature of the Levi-Civita connection is    
encoded in the Cartan curvature two-forms  
$$ \vec R = d \, \vec \Omega - 
{\textstyle \frac 1 2}\, \vec \Omega \times \vec \Omega,  
\tag 4.2 $$
and the metric is flat if and only if  $\vec R = 0$.

On $U \subset \comp^3$ with coordinates $(x_1,x_2,x_3)$,  
we say that a triplet of one-forms $\vec \theta$ 
is a ``canonical coframe" if there 
is a triplet $\vec H$ of  functions such that 
$$ 
\vec H = \pmatrix h_1 \cr h_2 \cr h_3 \cr \endpmatrix, \qquad 
\vec \theta = \pmatrix \theta_1 \cr 
\theta_2 \cr \theta_3 \cr \endpmatrix
= \pmatrix h_1 \, dx_1 \cr
h_2 \, dx_2 \cr
 h_1 \, dx_3 \cr
\endpmatrix, \qquad 
\vec e = \pmatrix e_1 \cr e_2 \cr e_3 \cr \endpmatrix = 
\pmatrix h_1^{-1} \, \partial_1 \cr 
h_2^{-1} \, \partial_2 \cr
h_3^{-1} \, \partial_3 \cr 
\endpmatrix, 
$$
where $\vec e$ is the orthonormal frame dual to $\vec\theta$.  
A canonical coframe can be reconstructed from 
its ``canonical coefficients" $\vec H = i_I \, \vec \theta$.  
Obviously a metric $g$ admits an orthonormal  
canonical coframe if and only if $g$ diagonal, eq.(2.3).      
We define a homogeneity-$m$ 
 ``Frobenius coframe" to be  a canonical coframe $\vec \theta$ such that 
$$ g = \theta_1 \otimes \theta_1 + \theta_2 \otimes \theta_2 + 
\theta_3 \otimes \theta_3 
= h_1^2 \, dx_1 \otimes dx_1 + h_2^2 \, dx_2 \otimes dx_2 + 
h_3^2 \, dx_3 \otimes dx_3 $$
is a homogeneity-$m$ Frobenius metric. 

We  say that a connection form $\vec \Omega$ is 
``Egoroff" if there is a triplet $\vec F$ of functions 
such that
$$ 
\vec F =   \pmatrix f_1 \cr f_2 \cr f_3 \endpmatrix 
, \qquad 
\vec \Omega = \pmatrix \Omega_1 \cr \Omega_2 \cr \Omega_3 \cr \endpmatrix 
= \pmatrix 
f_1 \, (x_2 - x_3)^{-1} \, (dx_2 - dx_3) \cr
f_2 \, (x_3 - x_1)^{-1} \, (dx_3 - dx_1) \cr
f_3 \, (x_1 - x_2)^{-1} \, (dx_1 - dx_2) \cr
\endpmatrix . \tag 4.3
$$
An Egoroff connection form can be 
reconstructed from the ``Egoroff coefficients"  
$\vec F= i_E \, \vec \Omega$.  
This  
is a nonstandard definition, but in Lemma A.3 we'll show that 
the Levi-Civita connection form $\vec \Omega$ of an orthonormal canonical coframe 
is Egoroff if and only if the metric $g$ is Egoroff in the   
standard sense.    

The following proposition, which is 
proved in the Appendix, simultaneously characterizes three-dimensional 
Frobenius coframes and their Levi-Civita connection forms:     

\proclaim{Proposition 4.1} 
Let $\vec \theta$ be a canonical coframe.   
Then $\vec \theta$ is a homogeneity-$m$ Frobenius coframe with 
Levi-Civita connection  form $\vec \Omega$ if and only if the  
following four conditions hold:  
\roster
\item"(C1)" $\vec \Omega$  is an Egoroff connection form. 
\item"(C2)" $d \, \vec F - \vec \Omega \times \vec F = 0 $,  
where $\vec F = i_E \, \vec \Omega$. 
\item"(C3)" $ d \, \vec H - \vec \Omega \times \vec H = 0 $, where 
$\vec H = i_I \, \vec \theta$.  
\item"(C4)" $\vec F \times \vec H = \frac m 2 \, \vec H $. 
\endroster 
\endproclaim

\noindent The remainder of this section is based on corollaries of 
Proposition 4.1.  

The following differential equation for  $\vec F$ 
will be called 
the ``structural equation":   
$$\eqalign{
d\, f_1 - f_2 \, f_3 \, \left( {d x_3 - dx_1 \over x_3 - x_1} -
{d x_1 - d x_2 \over x_1 - x_2 } \right) &= 0, \cr
d\, f_2 - f_3 \, f_1 \, \left( {d x_1 - dx_2 \over x_1 - x_2} -
{d x_2 - d x_3 \over x_2 - x_3 } \right) &= 0, \cr
d\, f_3 - f_1 \, f_2 \, \left( {d x_2 - dx_3 \over x_2 - x_3} -
{d x_3 - d x_1 \over x_3 - x_1 } \right) &= 0 .  \cr
}\tag 4.4
$$
If $\vec\theta$ is a Frobenius coframe with Levi-Civita connection form 
$\vec \Omega$,   then 
conditions (C1) and (C2) of Proposition 4.1 are equivalent 
to the statement that the Egoroff coefficients $\vec F = i_E \, \vec \Omega$ 
solve the structural equation.  The following 
corollary of Proposition 4.1 is only slightly less obvious. 

\proclaim{Corollary 4.2} 
A  canonical coframe $\vec \theta$ is a 
nontrivial homogeneity-$0$ Frobenius coframe if and only if 
some nonzero constant multiple  
$c\, \vec H = i_I (c \, \vec \theta)$ of the canonical coefficients 
is a solution of the structural equation.   
\endproclaim

\demo{Proof} 
Suppose $\vec \theta$ is a nontrivial homogeneity-$0$ Frobenius  
coframe, and $\vec H = i_I \, \vec \theta$.  
Then the Levi-Civita connection form $\vec \Omega$ is nonzero, 
and $\vec F = i_E \, \vec \Omega$ is a nonzero solution of 
the structural equation by (C1) and (C2).  
Now (C4) implies $c \, \vec H = \vec F$ for some 
function $c$, but (C2) and (C3) imply $d \, c = 0$, so 
$c$ is a  constant.  

Conversely, suppose $c$ is a nonzero constant and that 
$c \, \vec H$ is a solution of the structural equation.  
Let $\vec \theta$ be the canonical coframe with canonical 
coefficients $\vec H$, and 
let $\vec \Omega$ be the Egoroff connection form with 
Egoroff coefficients $c \, \vec H$.   
Since $\vec F = i_E \, \vec \Omega = c \, \vec H$, this data 
satisfies  
conditions (C1)-(C4) of Proposition 4.1,  with $m=0$. 
Therefore $\vec \theta$ is a homogeniety-$0$ Frobenius coframe 
(and $\vec \Omega$ is the Levi-Civita connection form).  
 \qed
\enddemo

It is now an easy matter to prove Proposition 2.2.  
A simple computation shows that 
 structural equation eq.(4.4) is equivalent to  
$$ 
d\, f_1  - {f_2 \, f_3 \over t} \, dt
 = 0, \qquad
d \, f_2  - {f_3 \, f_1 \over 1-t} \, dt
= 0, \qquad
d\, f_3  - {f_1 \, f_2 \over t\,(t-1)} \, dt = 0 , 
\tag 4.5
$$
where $t$ is the cross-ratio eq.(2.4).   
Dubrovin \cite{D; eq.(3.113)} had obtained eq.(4.5) from 
a Hamiltonian approach to Frobenius manifolds,  see also \cite{Hi1}.   
The same equation, or more precisely 
its reduction by the $B$-symmetry to an ODE, appears in the work of 
Tod \cite{T,Hi2} on 
Riemannian metrics with self-dual curvature in (real) dimension four.  
Proposition 2.2 is proved by rewriting  eq.(4.5) as 
$$ 
t \, d\, f_1^2   = 2 \, {f_1 \, f_2 \, f_3 } \, dt, 
\qquad
(1-t) \, d \, f_2^2  = 2 \, {f_1 \, f_3 \, f_1 } \, dt,
 \qquad
t\,(t-1) \, d\, f_3^2  = 2 \, {f_1 \, f_2  \,f_3} \, dt = 0  
$$
and appealing to Corollary 4.2 to write $c^2 \, g_{ii} = f_i^2$.   
The trace eq.(2.6) is expressed in terms of $\vec F$ by   
$$ k = - {\frac {c^2} 2} \, (g_{11} + g_{22} + g_{33}) =
- {\textstyle\frac 1 2} \, (f_1^2 + f_2^2 + f_3^2) =
- {\textstyle\frac 1 2} \, \vec F \cdot \vec F. 
$$
For any solution $\vec F$ of the structural equation, it 
is also clear from the orthogonality of the connection that 
the trace  $k = -\frac 1 2 \,\vec F \cdot \vec F$ is constant. 

The trace is related to Dubrovin's  $\mu$ by 
$k = \frac 1 2 \, \mu^2$, compare  
 \cite{D; eq.(3.114)}.   The three-dimensional Frobenius manifolds 
associated to the Coxeter groups $A_3$, $B_3$, and $H_3$ have   
trace  $k = \frac 1 {32}$, $\frac 2 {25}$, and 
$\frac 1 {18}$ respectively \cite{D; App. E}.  
These values do not appear on the list $k = \frac 1 2 (n + \frac 1 2)^2$, 
$n = 0,1,2,\dots$  
of Theorem 3.1.

Our final topic is the basic classification theory  
of three-dimensional Frobenius coframes.  
An equivalence class $\vec{ [ \theta ] }$ of 
coframes under the equivalence relation of constant rescaling  
 will be called a  ``projective coframe".   
It is evident from the Cartan torsion equation that 
the Levi-Civita connection form $\vec \Omega$ of a coframe $\vec \theta$ 
depends only on the projective coframe $\vec{ [ \theta ] }$, as 
does $\vec F = i_E \, \vec \Omega$.   
This defines a function   
${\Cal S}: \vec {[\theta ]} \mapsto \vec F$ mapping projective 
coframes to triplets of functions.  
It is clear that ${\Cal S}$ maps  Frobenius projective coframes to  
solutions of the structural equation. 

\proclaim{Corollary 4.3} 
\roster
\item
The map ${\Cal S}$ restricts to a 
bijection from nontrivial  homogeneity $m=0$ Frobenius projective  coframes to 
nonzero  solutions of the structural equation.   
\item For $m \ne 0$,   
the map ${\Cal S}$ restricts to a 
bijection from  
homogeneity-$m$ Frobenius projective coframes to 
solutions of the structural equation of trace $k=m^2/8$.  
\endroster
\endproclaim

\demo{Proof}  
\roster
\item
This follows immediately from Corollary 4.2.  
\item 
We first show that for nonzero $m$, ${\Cal S}$  
maps  homogeneity-$m$ Frobenius projective coframes to solutions 
of the structural equation of trace $k=m^2/8$. 
Suppose $\vec \theta$ is a homogeneity-$m$ Frobenius coframe.   
By  (C4),  $\vec H$ is pointwise 
an eigenvector with eigenvalue $\lambda = \frac m 2$ 
of the linear map \nobreak{$M : \vec H \mapsto \vec F \times \vec H$}.  The 
identity
$\vec F \times (\vec F \times \vec H) =
 (\vec F \cdot \vec H) \, \vec F - (\vec F \cdot \vec F) \, \vec H$ gives 
the characteristic equation  
$ M \, ( M^2 - 2 \, k ) = 0 $.  Now the 
 eigenvalue $\lambda = \frac m 2$ is nonzero, so 
$ \lambda \,( \lambda^2 - 2 \, k )=0$ implies  $k = m^2/8$.  

We next establish injectivity.  Suppose 
$\vec \theta$ and $\vec \theta '$ are 
homogeneity-$m$ Frobenius coframes with Levi-Civita 
connection forms $\vec \Omega$ and $\vec \Omega'$ 
respectively, and suppose that  
${\Cal S}(\vec{ [\theta]}) = {\Cal S}(\vec{ [\theta']})$.  
The Egoroff property (C1) of the connection forms then 
implies $\vec \Omega = \vec \Omega'$.  
By (C4), the canonical coefficients  
$\vec H = i_I \vec \theta$ and $\vec H'= i_I \vec \theta'$ 
are  pointwise eigenvectors of $M$ with eigenvalue $\frac m 2$.  
 $M$ has three distinct 
eigenvalues since $k = m^2/8$ is nonzero, so  $\vec H ' = c \, \vec H$ for 
some scalar function $c$, and (C3)  implies $d \, c = 0$.  
Since $\vec \theta' = c \, \vec \theta$ 
for a constant $c$,  we have shown that  
$\vec{ [\theta]} = \vec{ [\theta']}$.  

We finally establish surjectivity.  
 Let $\vec F$ be a solution of the structural equation
 of nonzero trace $k = m^2/8$.  
Let $\vec \Omega$ be the Egoroff connection form with 
Egoroff coefficients $\vec F$.  Then the Cartan curvature form eq.(4.2) 
vanishes, see eq.(A.1) below, so the connection is flat.   
Let $\vec F_p$ be the value of $\vec F$ at at a point $p \in U$, 
and  choose a nonzero 
$\vec H_p$ such that $\vec F_p \times \vec H_p = \frac m 2 \, \vec H_p$. 
Assuming that $U$ is simply connected, 
parallel transport with the flat connection generates the unique 
$\vec H$ that satisfies condition (C3) and has the value 
$\vec H_p$ at $p$.   
Furthermore $\vec H$ satisfies (C4) on $U$, because   
$\vec F \times \vec H - \frac m 2 \, \vec H$ is covariantly constant 
by (C2) and (C3), and  $\vec F_p \times \vec H_p - \frac m 2 \, \vec H_p =0$. 
This data satisfies conditions (C1)-(C4) of Proposition 4.1, so the 
canonical coframe $\vec \theta$ with canonical coefficients $\vec H$ 
is a homogeneity-$m$ Frobenius coframe such that  
${\Cal S}(\vec{ [\theta]}) = \vec F$.   
\qed
\endroster
\enddemo

\noindent 
To summarize,  a trace-$k$ solution $\vec F$  of the structural equation 
generates a homogeneity-$m$ Frobenius coframe for each 
distinct root $m$ of $m \, (m^2 - 8 \, k) = 0$.  
The explicit construction of the homogeneity-$0$ Frobenius coframe from 
$\vec F$ is trivial by Corollary 4.2.  The explicit construction of 
the nonzero homogeneity Frobenius coframes is somewhat more 
complicated, see Hitchin \cite{Hi} for details.

\head Appendix:~Proof of Proposition 4.1 \endhead

To prove Proposition 4.1, we need  to  
establish that the conditions 
 (C1)-(C4) together with the Cartan torsion equation 
are equivalent to the conditions (M1)-(M3).  

\proclaim{Lemma A.1} 
If (C1)-(C4) hold, then $d \, \vec \theta - \vec \Omega \times \vec \theta = 0$. 
\endproclaim
  
\demo{Proof} 
The proof only requires conditions (C1) and (C3).  
The first component of the torsion of a canonical coframe is 
$$ \eqalign{d \, \theta_1 - \Omega_2 \, \theta_3 + \Omega_3 \, \theta_2 &= 
 d( h_1 \, dx_1) - \Omega_2 \wedge (h_3 \, dx_3) + \Omega_3 \wedge( h_2 \, dx_2).  \cr
}
$$
The Egoroff condition (C1) is equivalent to 
$$ \Omega_1 \wedge (dx_2 - dx_3) = 0, \quad
\Omega_2 \wedge (dx_3 - dx_1) = 0, \quad
\Omega_3 \wedge (dx_1 - dx_2) = 0,  
$$
and we have 
$$ \eqalign{d \, \theta_1 - \Omega_2 \, \theta_3 + \Omega_3 \, \theta_2  
 &= (d \,  h_1 - \Omega_2 \,  h_3  + 
\Omega_3 \, h_2 ) \wedge dx_1,  \cr
}
$$
which vanishes by (C3).  The other components of the torsion 
vanish analogously.  
\qed
\enddemo

\noindent We may therefore assume 
$d \, \vec \theta - \vec \Omega \times \vec \theta = 0$, and 
establish ((C1)-(C4))$\Leftrightarrow$((M1)-(M3)) under this assumption.  
We will break this up into a number of separate steps:  
\roster
\item  (C3)$\Leftrightarrow$(M3)
\item  ((C3) and (C4))$\Leftrightarrow$((M3) and (M2))
\item  ((C1) and (C2))$\Rightarrow$(M1)
\item  (C2)$\Leftarrow$((M1) and (M2))
\item  (C3)$\Rightarrow$(C1)
\endroster

\demo{Proof of (C3)$\Leftrightarrow$(M3)}
 $I = h_1 \, e_1 + h_2 \, e_2 + h_3 \, e_3$, so the covariant constancy 
of $I$ is equivalent to $\vec \nabla \vec H = d \, \vec H - \vec \Omega 
\times \vec H = 0$.  
\qed\enddemo

\demo{Proof of ((C3) and (C4))$\Leftrightarrow$((M3) and (M2))} 
(M2) is equivalent to $\lie_E \vec H = {\textstyle \frac m 2} \, \vec H$. 
If (C3) or equivalently (M3) holds, then 
$$ 0 = i_E (d \vec H - \vec \Omega \times \vec H ) = 
\lie_E \vec H - \vec F \times \vec H, $$
and then (M2) is equivalent to  (C4).  
\qed\enddemo

\demo{Proof of ((C1) and (C2))$\Rightarrow$(M1)}
Substitute $\Omega_i$ from eq.(4.3) into the first component of 
the Cartan curvature eqation eq.(4.2) and use 
eq.(4.4); 
$$ \eqalign{
R_1 &= d \, \Omega_1 - \Omega_2 \wedge \Omega_3 \cr
&=  d \,  f_1 \wedge \left({dx_2 - dx_3  \over x_2 - x_3}\right) -
f_2 \, f_3 \, \left( { dx_3 - dx_1 \over x_3 - x_1} \right)
 \wedge \left( {dx_1 - dx_2 \over x_1 - x_2 }\right)  \cr
&= - f_2 \, f_3 \, 
{ i_E \, \left( (dx_2 - dx_3) \wedge (dx_3 - dx_1)\wedge(dx_1 - dx_2)  \right)
\over (x_2 - x_3)(x_3 -x_1)(x_1 - x_2) }, \cr
} \tag A.1
$$
which vanishes because $(dx_2 - dx_3) \wedge (dx_3 - dx_1)\wedge(dx_1 - dx_2)=0$. 
The other components of the curvature form $\vec R$ vanish similarly. 
\qed\enddemo

Recall that a vector field $B$ is said to be a ``conformal Killing vector" for 
the Riemannian metric $g$ if $\lie_B \, g = r \, g$ for some 
constant $r$.  The following standard lemma holds in any dimension.   

\proclaim{Lemma A.2} 
If $B$ is a conformal Killing vector for a flat Riemannian metric $g$, 
then the tangent bundle endomorphism $\nabla B$ is 
covariantly constant.    
\endproclaim

\demo{Proof} A flat Riemannian metric locally admits 
``flat coordinates" $\{ t_i \}$ such that 
$\{ \epsilon_i = \frac \partial {\partial t_i} \}$ 
is an orthonormal frame.  The basis vector fields are covariantly constant, 
$\nabla \epsilon_i = 0$, and the 
covariant derivatives along the basis vectors commute, 
$\nabla_{\epsilon_i} \, \nabla_{\epsilon_j} = 
\nabla_{\epsilon_j} \, \nabla_{\epsilon_i}$.  
The metric equals $g = \sum_j dt_j \otimes dt_j$, where 
$\{ dt_i \}$ is the dual coframe.

Writing   
$ B = \sum_k b_k \, \epsilon_k$, we have  
$\lie_B \, dt_j = d ( \lie_B \, t_j) 
= d \, b_j 
= \sum_k (\nabla_{\epsilon_k} b_j) \, dt_k 
$, and   
$$\eqalign{ \lie_B \, g &= \sum_j \, \lie_B \left( dt_j \otimes dt_j\right) 
= \sum_{j,k} (\nabla_{\epsilon_j} b_k + \nabla_{\epsilon_k} b_j) 
\, dt_j \otimes dt_k.  \cr}
$$
Since $r$ is constant and $g$ is covariantly constant, 
$\lie_B \, g = r \, g$ implies 
$$ \eqalign{
0 &= \nabla_{\epsilon_i} ( \lie_B \, g ) = \sum_{j,k}
(\nabla_{\epsilon_i} \nabla_{\epsilon_j} b_k + 
\nabla_{\epsilon_k} \nabla_{\epsilon_i}  b_j) \, dt_j \otimes dt_k. \cr
}
$$
After several permutations of the indices,  
$$ \nabla_{\epsilon_i} \, \nabla_{\epsilon_j} \, b_k = 
- \nabla_{\epsilon_k} \, \nabla_{\epsilon_i} \, b_j =
 \nabla_{\epsilon_j} \, \nabla_{\epsilon_k} \, b_i =
-\nabla_{\epsilon_i} \, \nabla_{\epsilon_j} \, b_k ,  $$
which implies $\nabla_{\epsilon_i} \, \nabla_{\epsilon_j} \, b_k = 0$ 
for any $i,j,k$.  We conclude that $\nabla B$ is covariantly constant, as  
$$ \nabla_{\epsilon_i} (\nabla B) = 
\sum_{j,k} (\nabla_{\epsilon_i} \, \nabla_{\epsilon_j} \, b_k ) \, dt_j 
= 0 . \qed
$$
\enddemo

\demo{Proof of (C2)$\Leftarrow$((M1) and (M2))} 
Assume (M1) and (M2) hold.  
Since $\lie_E \, dx_i = dx_i$,  (M2) implies   
$ \lie_E \, g = \left(m+2 \right) \, g $. 
From  (M1) and Lemma A.2 we conclude 
that $\nabla E$ is covariantly constant.  
We will now show that $\nabla E$ is covariantly constant only if (C2) holds.  
  
We  compute $ \nabla_Z (\nabla E )$ 
for an arbitrary vector field $Z$.  From the Leibniz property 
$$\eqalign{ \nabla_Z \langle \nabla E , Y \rangle  &= 
 \langle \nabla_Z (\nabla E), Y \rangle + \langle \nabla E , \nabla_Z Y \rangle \cr
} $$
conclude that 
$$\eqalign{  \langle \nabla_Z (\nabla E), Y \rangle &=
   \nabla_Z \nabla_Y E -  \nabla_{(\nabla_Z Y)} E.  \cr
} \tag A.2 $$
We compute $\nabla_V E$ for an arbitrary vector 
field $V$.  The vanishing torsion of the Levi-Civita 
connection gives 
$$\eqalign{ \nabla_V E & = \nabla_E V -\lie_E V . \cr
}\tag A.3
$$
Now $e_i = h_i^{-1} \, \partial_i$, so using  
the Leibniz property of the Lie derivative  and 
the homogeneity property $\lie_E \, h_i^{-1} = - \frac m 2 \, h_i^{-1}$, we have  
 $\lie_E \, e_i = -\beta \, e_i$,  
where $\beta = \left( \frac m 2 + 1 \right)$, so
$$\eqalign{
\lie_E V &= \lie_E (V_1 \, e_1 + V_2 \, e_2 + V_3 \, e_3) \cr
&= 
(\lie_E V_1 - \beta\,  V_1)\,  e_1 + (\lie_E V_2 - \beta\,  V_2)\,  e_2
+(\lie_E V_3 - \beta\,  V_3)\,  e_3.  \cr
}
$$
Now eq.(A.3) written in in terms of components is  
$$
\eqalign{\vec\nabla_V \vec E 
&= 
(\lie_E \vec V - (i_E \vec \Omega) \times \vec V)  - 
(\lie_E \vec V -\beta \, \vec V )\cr
&= - \vec F \times \vec V + \beta\,  \vec V \cr
}\tag A.4 
$$
Setting $V = Y$ in eq.(A.4) and applying $\vec \nabla_Z$ gives 
$$\vec \nabla_Z \vec \nabla_Y \vec E  = - 
(\vec \nabla_Z \vec F) \times \vec Y - \vec F \times (\vec \nabla_Z \vec Y) 
+ \beta\, \vec \nabla_Z \vec Y , $$
while setting $V = \nabla_Z Y$ in eq.(A.4) gives 
$$\vec \nabla_{(\nabla_Z Y)} \vec E  = - 
\vec F \times (\vec \nabla_Z \vec Y) + 
\beta \,  \vec \nabla_Z \vec Y,  $$
and from eq.(A.2) we conclude that
$  \langle \nabla_Z (\nabla E), Y \rangle = 0$ if and only if 
$(\vec \nabla_Z \vec F) \times \vec Y = 0$.   Since $Y$ and $Z$ 
are arbitrary,  $\nabla E$ is covariantly constant 
if and only if $0 = \vec \nabla \vec F  = d \, \vec F - \vec \Omega \times 
\vec F$. This is just the condition (C2).  
\qed 
\enddemo

The following lemma will complete the proof of Proposition 4.1, 
and also establish the equivalence of our definition  of 
the Egoroff condition, eq.(4.3), with the usual definition, (C1a) below.

\proclaim{Lemma A.3} Let $\vec \theta$ be a canonical coframe 
with canonical coefficients $\vec H = i_I \, \vec \theta$, and 
let $\vec \Omega$ be its Levi-Civita connection form.  
Then the following  are equivalent.
\roster 
\item"(C1)" 
$\vec \Omega$ is an Egoroff connection form.  
\item"(C1a)" The one-form
$ \vec H \cdot \vec \theta =   
g_{11} \, dx_1 + g_{22} \, dx_2 + g_{33} \, dx_3 $ is closed. 
\item"(C1b)" $(d \vec H - \vec \Omega \times \vec H ) \cdot \vec \theta = 0. $ 
\item"(C1c)"
$i_I \, \vec \Omega = 0$.  

\endroster
\endproclaim

\demo{Proof} 
The equivalence (C1a)$\Leftrightarrow$(C1b) is immediate from the derivation 
property of $d$ and the 
triple-product identity:  
$$ \eqalign{
d\, (\vec H \cdot \vec \theta) 
&= (d \vec H) \cdot \vec \theta + \vec H \cdot (d \vec \theta) 
= (d \vec H) \cdot \vec \theta + \vec H \cdot 
(\vec \Omega \times \vec \theta) 
= (d \vec H - \vec \Omega \times \vec H) \cdot \vec \theta . \cr
}
$$
To prove (C1b)$\Leftrightarrow$(C1c), start with 
$$ \eqalign{ 
d \vec H - \vec \Omega \times \vec H 
&= d \, i_I \vec \theta - \vec \Omega \times (i_I \vec \theta) \cr
&= (\lie_I \, \vec \theta - i_I \, d \, \vec \theta ) +
i_I \, (\vec \Omega \times \vec \theta) - (i_I \vec \Omega) \times \vec \theta
\cr
&= \lie_I \, \vec \theta -  (i_I \vec \Omega) \times \vec \theta
-i_I \left( \, d \, \vec \theta  -\vec \Omega \times \vec \theta \right)
\cr
&= \lie_I \, \vec \theta -  (i_I \vec \Omega) \times \vec \theta, 
\cr
}
$$
which implies  
$$\eqalign{ (d \vec H - \vec \Omega \times \vec H) \cdot \vec \theta 
&=(\lie_I \, \vec \theta -  (i_I \vec \Omega) \times \vec \theta) \cdot 
\vec \theta = 
(\lie_I \, \vec \theta) \cdot \vec \theta
 -  (i_I \vec \Omega) \cdot (\vec \theta \times
\vec \theta) .\cr
} 
$$
Now  
$$(\lie_I \, \vec \theta) \cdot \vec \theta = 
 (\lie_I \, h_1) \, h_1 \, dx_1 \wedge dx_1 
+
(\lie_I \, h_2) \, h_2 \, dx_2 \wedge dx_2+
(\lie_I \, h_2) \, h_2 \, dx_2 \wedge dx_2 = 0,  $$
and 
$$ (i_I \vec \Omega) \cdot (\vec \theta \times
\vec \theta) = 2 \, (i_I  \Omega_1) \, \theta_2 \wedge \theta_3 +
2 \, (i_I  \Omega_2) \, \theta_3 \wedge \theta_1  +
2 \, (i_I  \Omega_3) \, \theta_1 \wedge \theta_2 , 
$$
so $(d \vec H - \vec \Omega \times \vec H) \cdot \vec \theta$ vanishes 
if and only if each component of $i_I  \vec \Omega$ vanishes. 

Finally we prove (C1c)$\Leftrightarrow$(C1). 
Any coframe $\vec \theta$ and dual frame $\vec e$  tautologically satisfy 
$i_{e_i} \, \theta_j = \delta_{ij}$.  A canonical coframe 
has the additional property $i_{e_2} \, i_{e_3} \, d \, \theta_1 = 0$, 
which follows from 
$d \, \theta_1 = h_1^{-1} \, d \, h_1 \wedge \theta_1$.  
The first component of the Cartan torsion equation gives 
$$ \eqalign{
0 &= i_{e_2} \, i_{e_3} \, 
(d \, \theta_1 - \Omega_2 \wedge \theta_3 + \Omega_3 \wedge \theta_2) =
 i_{e_2} \, \Omega_2 + i_{e_3} \, \Omega_3,  \cr
}
$$
which together with the 
other components  $i_{e_3} \, \Omega_3 + i_{e_1} \, \Omega_1= 0$ 
and $i_{e_1} \, \Omega_1 + i_{e_2} \, \Omega_2=0$  implies 
$i_{e_i} \, \Omega_i = 0$, or equivalently $i_{\partial_i} \, \Omega_i = 0$.  
So the connection form of a canonical coframe satisfies 
$$ 
\vec \Omega = 
\pmatrix 
i_{\partial_2} \, \Omega_1 \, dx_2 + i_{\partial_3} \, \Omega_1 \, dx_3\cr
i_{\partial_3} \, \Omega_2 \, dx_3 + i_{\partial_1} \, \Omega_2 \, dx_1\cr
i_{\partial_1} \, \Omega_3 \, dx_1 + i_{\partial_2} \, \Omega_3 \, dx_2\cr
\endpmatrix, \qquad
i_I \, \vec \Omega = 
\pmatrix 
i_{\partial_2} \, \Omega_1 + i_{\partial_3} \, \Omega_1 \cr
i_{\partial_3} \, \Omega_2 + i_{\partial_1} \, \Omega_2 \cr
i_{\partial_1} \, \Omega_3 + i_{\partial_2} \, \Omega_3 \cr
\endpmatrix. 
$$
Now (C1c) is equivalent to 
$$i_{\partial_3} \, \Omega_1 = -  i_{\partial_2} \, \Omega_1, \qquad
i_{\partial_1} \, \Omega_2 = -  i_{\partial_3} \, \Omega_2, \qquad
i_{\partial_2} \, \Omega_3 = -  i_{\partial_1} \, \Omega_3, \qquad
$$
which is equivalent to (C1),  compare with eq.(4.3):  
$$  \vec \Omega = 
\pmatrix 
f_1 \, (x_2 - x_3)^{-1} \, (dx_2 - dx_3) \cr
f_2 \, (x_3 - x_1)^{-1} \, (dx_3 - dx_1) \cr
f_3 \, (x_1 - x_2)^{-1} \, (dx_1 - dx_2) \cr
\endpmatrix , \qquad
\vec F = \pmatrix f_1 \cr f_2 \cr f_3 \endpmatrix =
\pmatrix 
 (x_2 -x_3)\, i_{\partial_2} \, \Omega_1 \cr
 (x_3 -x_1) \,i_{\partial_3} \, \Omega_2 \cr
 (x_1 -x_2) \, i_{\partial_1} \, \Omega_3 \cr
\endpmatrix .  \qed
$$
\enddemo 

\demo{Proof of (C3)$\Rightarrow$(C1)}  Obviously 
(C3)$\Rightarrow$(C1b), and (C1b)$\Leftrightarrow$(C1) by Lemma A.3. 
\qed\enddemo 

\noindent This completes the proof of Proposition 4.1.

\refstyle{A}

\enddocument

\widestnumber\key{ADHM}


\Refs

\ref\key  ADHM
\by M.F. Atiyah, V.G. Drinfeld, N.J. Hitchin, and Yu.I. Manin
\paper Construction of Instantons
\jour Phys. Lett. \vol 65A \page 185 \yr 1978 \endref

\ref\key  At 
\by M.F. Atiyah 
\book Geometry of Yang-Mills Fields
\publ Lezioni Fermiane, Accademia Nazionale dei Lince \& Scuola 
Normale
Superiore \publaddr Pisa \yr 1979 \endref

\ref\key Au
\by M. Audin
\paper An introduction to Frobenius manifolds, 
moduli spaces of stable maps and quantum cohomology
\paperinfo Preprint
\endref



\ref\key BS
\by G. Bor and J. Segert
\paper Symmetric instantons and the ADHM construction
\jour Commun. Math. Phys \vol 183 \pages 183--203 \yr 1997 \endref

\ref\key CDD
\by Y. Choquet-Bruhat, C. DeWitt-Morette, and M. Dillard-Bleick
\book Analysis, Manifolds and Physics, revised edition
\publ North-Holland \publaddr Amsterdam \yr 1982
\endref

\ref\key D
\by B. Dubrovin
\paper Geometry of 2D Topological Field Theories
\inbook Lecture Notes in Mathematics \vol 1620 
\pages 120--348 \publ Springer-Verlag \yr 1996 
\endref

\ref\key F
\by R. Fuchs
\paper \"Uber lineare homogene Differentialgleichungen zweiter 
Ordnung mit drei im Endlichen gelegenen wesentlich wingul\"aren Stellen
\jour Math. Ann. \vol 63 \pages 301--321 \yr 1907 \endref



\ref\key G
\by A.B. Givental
\paper Equivariant Gromov-Witten invariants
\jour Internat. Math. Res. Notices 
\vol 13 \yr 1996 \pages 613--663 
\endref 


\ref\key Hi1
\by N.J. Hitchin
\paper Frobenius manifolds (Notes by David Calderbank)
\paperinfo Preprint  
\yr 1996 
\endref

\ref\key Hi2
\by N.J. Hitchin
\paper Twistor spaces, Einsten metrics and isomonodromic deformations
\jour J. Diff. Geom. \vol 42 \pages 30--112 \yr 1995 
\endref

\ref\key Hi3
\by N.J. Hitchin
\paper Poncelet polygons and the Painlev\'e equations
\inbook Geometry and Analysis
\pages 151--185
\bookinfo Tata Institute of Fundamental Research, Bombay 
\endref



\ref\key JM
\by M. Jimbo and T. Miwa
\paper Monodromy preserving deformations of linear ordinary differential 
equations with rational coefficients. II
\jour Physica \vol 2D \pages 407--448 \yr 1981 
\endref

\ref\key Ml
\by B. Malgrange
\paper Sur les deformations isomonodromiques I. Singularit\'es r\'eguli\`eres
\inbook Math\'ematique et Physique, S\'eminaire de l'Ecole Normale Sup\'erieure 
1979-1982 \pages 401--426 
\bookinfo Progress in Mathematics \vol 37
\publ Birkh\"auser \yr 1983 
\endref

\ref\key Mn1 
\by Yu.I. Manin
\paper Frobenius manifolds, quantum cohomology, and moduli spaces
\paperinfo Preprint MPI 96-113
\endref


\ref\key Mn2
\by Yu.I. Manin
\paper Sixth Painlev\'e equation, universal elliptic curve, and mirror of $P^2$
\paperinfo Preprint, al-geom/9605010 \yr 1996 
\endref

\ref\key RT
\by Y. Ruan and G. Tian
\paper A mathematical theory of quantum cohomology
\jour J. Diff. Geom 
\vol42 \yr 1995 \pages 269--367
\endref

\ref\key Sa
\by C. Sabbah
\paper Frobenius manifolds: isomonodromic deformations and 
infinitesimal period mappings
\paperinfo Preprint 
\endref

\ref\key Se
\by J. Segert
\paper Painlev\'e solutions from equivariant 
holomorphic bundles
\paperinfo Preprint
\yr 1996 \endref  

 \ref\key T
\by K.P. Tod
\paper Self-dual Einstein metrics from the Painlev\'e VI equation
\jour Phys. Lett A \vol 190 \yr 1994 \pages 221-224 \endref
 

\endRefs
\bye